\documentclass[11pt]{amsart}
\usepackage{graphicx}
\usepackage{amssymb, amsmath}
\theoremstyle{definition}

\theoremstyle{remark}

\numberwithin{equation}{section}

\begin{document}
\title{A note on cellular automata}
\author{M. Shahryari}
\address{M. Shahryari\\
Department of Pure Mathematics,  Faculty of Mathematical
Sciences, University of Tabriz, Tabriz, Iran }
\email{mshahryari@tabrizu.ac.ir}
\date{\today}

\begin{abstract}
For an arbitrary group $G$ and arbitrary set $A$, we define a monoid structure on the set of all uniformly continuous functions $A^G\to A$ and then we show that it is naturally isomorphic  to the monoid  of cellular automata $\mathrm{CA}(G, A)$. This gives a new equivalent definition of a cellular automaton over the group $G$ with alphabet set $A$. We use this new interpretation to give a simple proof of the theorem of Curtis-Hedlund.
\end{abstract}

\maketitle


Let $G$ be a group and $A$ be a non-empty  set. The set $A^G$ consists of all configurations $x:G\to A$. We suppose that the sets $A^G$ and $A$ are equipped with the pro-discrete and discrete uniform structures, respectively (see \cite{Bour} for basic definitions). We also suppose that the group $G$ acts on $A^G$ by shift: $(g\cdot x)(h)=x(g^{-1}h)$, for all $g, h\in G$ and $x\in A^G$. Recall that a cellular automaton over $G$ with alphabet set $A$ is a map $T:A^G\to A^G$ such that there exists a finite subset $S\subseteq G$ and a function $\mu:A^S\to A$ with the following property: for all $x\in A^G$ and all $g\in G$, we have
$$
T(x)(g)=\mu((g^{-1}\cdot x)_{|_S}),
$$
where ${|_S}$ denotes the restriction. Any such a set $S$ is a called a memory set and $\mu$ is called a local defining function. Let $\mathrm{CA}(G, A)$ be the set of all such cellular automata. This set is a monoid with the ordinary composition of mappings. The reader may see \cite{CS-C} for detailed discussion.

Let $\mathcal{U}(A^G, A)$ be the set of all uniformly continuous functions from $A^G$ to $A$. Define a binary operation $\ast$ on this set by
$$
(f_1\ast f_2)(x)=f_1((f_2(g^{-1}\cdot x))_{g\in G}).
$$
We soon will see that $\mathcal{U}(A^G, A)$ is a monoid. The identity of this monoid is the projection map $p_1:A^G\to A$ defined by $p_1(x)=x(1)$, where $1$ is the identity of $G$.   For any $T\in \mathrm{CA}(G, A)$, we define a new map $f_T:A^G\to A$ by
$$
f_T(x)=T(x)(1).
$$

{\bf Theorem A.} {\it
The map $T\mapsto f_T$ is an isomorphism between the monoids $\mathrm{CA}(G,A)$ and $\mathcal{U}(A^G, A)$.
}

\begin{proof}
We first show that $f_T$ is uniformly continuous. Suppose $S$ is a memory set for $T$ and put
$$
W_S=\{ (x,y)\in A^G\times A^G:\ x_{|_S}=y_{|_S}\}.
$$
Let $\mathcal{U}$ denote the pro-discrete uniform structure over $A^G$. We know that $W_S\in \mathcal{U}$. Consider the map $f_T\times f_T:A^G\times A^G\to A\times A$ defined by
$$
(f_T\times f_T)(x,y)=(f_T(x), f_T(y)).
$$
We have the implication
$$
x_{|_S}=y_{|_S}\Rightarrow T(x)(1)=T(y)(1),
$$
and this means that $W_S\subseteq (f_T\times f_T)^{-1}(\Delta_A)$, where $\Delta_A$ is the diagonal of $A\times A$. This shows that
$$
(f_T\times f_T)^{-1}(\Delta_A)\in \mathcal{U},
$$
and hence $f_T$ is uniformly continuous. Now, for arbitrary automata $T_1$ and $T_2$, we have
\begin{eqnarray*}
f_{T_1\circ T_2}(x)&=& (T_1\circ T_2)(x)(1)\\
                   &=&T_1(T_2(x))(1)\\
                   &=&f_{T_1}(T_2(x))\\
                   &=&f_{T_1}((T_2(x)(g))_{g\in G})\\
                   &=&f_{T_1}((g^{-1}\cdot T_2(x)(1))_{g\in G})\\
                   &=&f_{T_1}((T_2(g^{-1}\cdot x)(1))_{g\in G})\\
                   &=&f_{T_1}((f_{T_2}(g^{-1}\cdot x))_{g\in G})\\
                   &=& (f_{T_1}\ast f_{T_2})(x).
\end{eqnarray*}
This shows that the map $T\mapsto f_T$ is a homomorphism. Note that if $f_{T_1}=f_{T_2}$, then for any $x\in A^G$, we have $T_1(x)(1)=T_2(x)(1)$, and since $T_1$ and $T_2$ are $G$-equivariant, so $T_1=T_2$, proving that the map is injective.

Now, suppose that $f\in \mathcal{U}(A^G, A)$. Define a map $T:A^G\to A^G$ by $T(x)(g)=f(g^{-1}\cdot x)$. First, note that $T$ is $G$-equivariant: for any $h\in G$, we have
\begin{eqnarray*}
T(h\cdot x)(g)&=&f(g^{-1}\cdot (h\cdot x))\\
              &=&f((h^{-1}g)^{-1}\cdot x)\\
              &=&T(x)(h^{-1}g)\\
              &=&(h\cdot T(x))(g).
\end{eqnarray*}
Since $f$ is uniformly continuous, we have $(f\times f)^{-1}(\Delta_A)\in \mathcal{U}$. On the other hand, we know that the set
$$
\{ W_{\Omega}:\ \Omega\subseteq G, |\Omega|<\infty\}
$$
is a basis for the  pro-discrete uniform structure over $A^G$. Hence there exists a finite subset $S\subseteq G$ such that $W_S\subseteq (f\times f)^{-1}(\Delta_A)$. In other words
$$
x_{|_S}=y_{|_S}\Rightarrow f(x)=f(y) \Rightarrow T(x)(1)=T(y)(1).
$$
This shows that $T$ is a cellular automaton with the memory set $S$. Clearly $f_T=f$ and this shows that the map $T\mapsto f_T$ is surjective. Therefore, we proved that $\mathcal{U}(A^G, A)$ is a monoid and it is isomorphic to $CA(G,A)$.
\end{proof}

As a result, we now have a very easy definition of a cellular automaton: any uniformly continuous map $f:A^G\to A$ is a cellular automaton! As an application, we reprove the theorem of Curtis and Hedlund (see \cite{CS-C}):\\

{\bf Corollary A.} {\it
Let $A$ be finite and $T:A^G\to A^G$ be continuous and $G$-equivariant. Then $T$ is a cellular automaton.
}

\begin{proof}
Define a map $f_T:A^G\to A$ by $f_T(x)=T(x)(1)$. Note that $f_T=p_1\circ T$, so it is continuous. Since $A^G$ is compact, so $f_T$ is uniformly continuous and hence $f_T\in \mathcal{U}(A^G, A)$. This shows that there exists a cellular automaton $T_0$ such that $f_T=f_{T_0}$. But since $T$ is $G$-equivariant, it can be easily seen that $T_0=T$, proving that $T$ is a cellular automaton.
\end{proof}

As another application, we prove the existence of the {\em minimal memory set} for a cellular automaton.\\

{\bf Corollary B.}
{\it Let $T$ be a cellular automaton and $S$ and $S^{\prime}$ be two memory sets for $T$. Then $S\cap S^{\prime}$ is also a memory set.
}

\begin{proof}
Let $f:A^G\to A$ be the corresponding uniformly continuous mapping and $V=(f\times f)^{-1}(\Delta_A)$. We know that a finite set $\Omega\subseteq G$ is a memory set for $T$ if and only if $W_{\Omega}\subseteq V$. Clearly $V\circ V\subseteq V$. Let $(x, y)\in W_{S\cap S^{\prime}}$ and choose $z\in A^G$ such that $(x,z)\in W_S$ and $(z, y)\in W_{S^{\prime}}$. This shows that $(x,y)\in W_S\circ W_{S^{\prime}}$. Hence
$$
W_{S\cap S^{\prime}}\subseteq  W_S\circ W_{S^{\prime}}\subseteq V\circ V\subseteq V.
$$
This proves that $S\cap S^{\prime}$ is a memory set.
\end{proof}

Finally, we must say that a similar statement is true for any arbitrary subshift $X\leq A^G$, after a small modification: Let $\mathcal{U}_0(X, A)$ be the set of all uniformly continuous functions $f:X\to A$, with the further property that
$$
(f(g^{-1}\cdot x))_{g\in G}\in X,
$$
for all $x\in X$. Then we can define the binary operation
$$
(f_1\ast f_2)(x)=f_1((f_2(g^{-1}\cdot x))_{g\in G}).
$$
on the set $\mathcal{U}_0(X, A)$ and it becomes a monoid again. We can prove then the next theorem.\\

{\bf Theorem B.} {\it There is a natural isomorphism between $\mathcal{U}_0(X, A)$ and the monoid of all cellular automata $X\to X$.}


\begin{thebibliography}{99}

\bibitem{Bour} N. Bourbaki. {\it Elements of Mathematics: General Topology (Part 1)}. English translation, Springer.

\bibitem{CS-C} T. Ceccherini-Silberstein, M. Coornaert. {\it Cellular Automata and Groups}. Springer, 2010.




\end{thebibliography}
\end{document}